\documentclass[11pt]{article}
\usepackage{amssymb, amsthm, amsmath}

\newtheorem{thm}{Theorem}[section]
\newtheorem{lem}[thm]{Lemma}

\renewcommand{\phi}{\varphi}

\newcommand{\N}{\mathbb{N}}
\newcommand{\Z}{\mathbb{Z}}

\newcommand{\Homeo}{\operatorname{Homeo}}
\newcommand{\id}{\operatorname{id}}

\newcommand{\xp}{(X,\phi)}
\newcommand{\yp}{(Y,\psi)}

\title{A short proof of affability \\
for certain Cantor minimal $\Z^2$-systems}
\author{MATUI Hiroki}
\date{}

\begin{document}
\maketitle

\begin{abstract}
We will show that any extension of a product of two Cantor minimal 
$\Z$-systems is affable in the sense of Giordano, Putnam and Skau. 
\end{abstract}

\section{Introduction}

In this paper, we would like to investigate the orbit structure 
of certain minimal dynamical systems on a Cantor set. 
Giordano, Putnam and Skau proved that 
equivalence relations arising from $\Z$-actions are 
orbit equivalent to AF equivalence relations in \cite{GPS1}. 
Moreover, they gave the classification for AF equivalence relations. 
In a recent paper \cite{GPS3}, 
they continued their investigations and 
showed that equivalence relations arising from $\Z^2$-actions 
are again orbit equivalent to AF equivalence relations 
under a hypothesis involving the existence of cocycles. 
An equivalence relation which is orbit equivalent to 
an AF equivalence relation is said to be affable. 
A crucial ingredient of their proof was the absorption theorem 
given in \cite{GPS2}. 
They needed, however, sufficiently many cocycles 
in order to construct an AF subequivalence relation 
to which the absorption theorem can be applied. 
The aim of this paper is to show that the existence of cocycles 
is not necessary for certain $\Z^2$-actions. 
We will give a short proof 
that the associated equivalence relations 
are orbit equivalent to AF equivalence relations, 
thus they are affable. 

We recall some terminology that we shall use. 
Let $X$ be a Cantor set 
and let ${\cal R}$ be an \'etale equivalence relation on $X$. 
We define the ${\cal R}$-equivalence class $[x]_{\cal R}$ 
of $x\in X$ by $[x]_{\cal R}=\{y\in X:(x,y)\in{\cal R}\}$. 
The equivalence relation ${\cal R}$ is said to be minimal, 
if $[x]_{\cal R}$ is dense in $X$ for every $x\in X$. 
Let $\phi:G\to\Homeo(X)$ be a free action of 
a countable discrete group $G$, that is, 
$\phi$ is a group homomorphism and 
$\phi^g(x)\neq x$ for all $x\in X$ and $g\in G\setminus\{e\}$, 
where $e$ means the identity element. 
We put 
\[ {\cal R}_\phi=\{(x,\phi^g(x))\in X\times X:x\in X, \ g\in G\}. \]
By transferring the product topology on $X\times G$ 
via the bijection $(x,g)\mapsto(x,\phi^g(x))$, 
we can topologize ${\cal R}_\phi$. 
It is easily verified that 
${\cal R}_\phi$ becomes an \'etale equivalence relation. 
We call $\xp$ a Cantor minimal $G$-system 
when ${\cal R}_\phi$ is minimal. 
In this paper, we deal with only Cantor minimal $\Z$-systems 
and Cantor minimal $\Z^2$-systems. 

Let $\xp$ and $\yp$ be two Cantor minimal $\Z^2$-systems. 
We say that $\pi:\yp\to\xp$ is a factor map 
when $\pi:Y\to X$ is a continuous map and 
$\pi\circ\psi^a=\phi^a\circ\pi$ for all $a\in\Z^2$. 
The system $\yp$ is called an extension of $\xp$. 
Our main theorem asserts that ${\cal R}_\psi$ is affable, 
if $\xp$ is conjugate to a product of two Cantor minimal $\Z$-systems. 
Suppose that $\xp$ is conjugate to the product of 
two Cantor minimal $\Z$-systems $(X_1,\phi_1)$ and $(X_2,\phi_2)$. 
From \cite[Theorem 2.3]{GPS1} we can see that 
${\cal R}_{\phi_1}$ and ${\cal R}_{\phi_2}$ are affable. 
Since a product of AF equivalence relations is also AF, 
it is easily checked that ${\cal R}_\phi$ is affable. 
But, it looks impossible to mimic this simple argument 
in the case of the extension $\yp$, 
because $\yp$ itself is not a product. 
We will instead construct a `nice' AF subequivalence relation 
of ${\cal R}_\phi$ and 
apply the absorption theorem to this relation.

\section{Products of Cantor minimal $\Z$-systems}

Throughout this section, 
let $B_i=(V_i,E_i)$ be simple properly ordered Bratteli diagrams 
for $i=1,2$. 
For each $i=1,2$, $V_i$ and $E_i$ can be written 
as a countable disjoint union of non-empty finite sets: 
\[ V_i=V_{i,0}\cup V_{i,1}\cup V_{i,2}\cup\dots
\text{ and }
E_i=E_{i,1}\cup E_{i,2}\cup E_{i,3}\cup\dots \]
with the source map $s:E_{i,n}\to V_{i,n-1}$ and 
the range map $r:E_{i,n}\to V_{i,n}$. 
Without loss of generality, we may assume that 
all two vertices in consecutive levels are 
connected by more than three edges. 
We write the infinite path space associated with $B_i$ 
by $X_i$ for each $i=1,2$. 
Let $p_i$ be the unique maximal infinite path of $X_i$ and 
let $\phi_i\in\Homeo(X_i)$ be 
the Bratteli-Vershik transformation on $X_i$ (see \cite{HPS}). 
It is well known that 
$(X_i,\phi_i)$ is a Cantor minimal $\Z$-system. 

Set $X=X_1\times X_2$. 
Let $\phi:\Z^2\to\Homeo(X)$ be the $\Z^2$-action on $X$ 
induced by $\phi_1\times\id$ and $\id\times\phi_2$. 
Then, $\xp$ is a Cantor minimal $\Z^2$-system. 

For each $n\in\N$, we put 
\[ {\cal R}_n=\{((x_1,x_2),(y_1,y_2))\in X\times X:
x_{i,m}=y_{i,m}\text{ for }i=1,2
\text{ and }m>n\}, \]
where $x_{i,m},y_{i,m}\in E_{i,m}$ mean the $m$-th coordinate of 
the infinite paths $x_i,y_i\in X_i$. 
It is not hard to see that 
${\cal R}_n$ is a compact open subequivalence relation 
of ${\cal R}_\phi$ 
with the relative topology from ${\cal R}_\phi$. 
Therefore 
\[ {\cal R}=\bigcup_{n\in\N}{\cal R}_n \]
is an AF subequivalence relation of ${\cal R}_\phi$. 
Note that ${\cal R}$ is minimal because $B_1$ and $B_2$ are simple. 

For $i,j=0,1$ and $n\in\N$, we define continuous functions 
$\lambda^{ij}_n:X\to\{0,1\}$ inductively as follows. 
Let $(x_1,x_2)\in X$. 
We denote the $n$-th coordinate of $x_i$ by $x_{i,n}\in E_{i.n}$. 
At first, put 
\[ \lambda^{00}_1(x_1,x_2)=\begin{cases}
1 & x_{1,1}\text{ is maximal} \\
0 & \text{otherwise,} \end{cases} \]
\[ \lambda^{01}_1(x_1,x_2)=\begin{cases}
1 & x_{2,1}\text{ is minimal} \\
0 & \text{otherwise,} \end{cases} \]
\[ \lambda^{11}_1(x_1,x_2)=\begin{cases}
1 & x_{1,1}\text{ is minimal} \\
0 & \text{otherwise} \end{cases} \]
and 
\[ \lambda^{10}_1(x_1,x_2)=\begin{cases}
1 & x_{2,1}\text{ is maximal} \\
0 & \text{otherwise.} \end{cases} \]
Then, for $n\geq2$, we define $\lambda^{ij}_n$ by 
\[ \lambda^{00}_n(x_1,x_2)=\begin{cases}
\lambda^{00}_{n-1}(x_1,x_2) & 
x_{1,n}\text{ is maximal and }x_{2,n}\text{ is maximal} \\
1 & x_{1,n}\text{ is maximal and }x_{2,n}\text{ is not maximal} \\
0 & \text{otherwise,} \end{cases} \]
\[ \lambda^{01}_n(x_1,x_2)=\begin{cases}
\lambda^{01}_{n-1}(x_1,x_2) & 
x_{1,n}\text{ is maximal and }x_{2,n}\text{ is minimal} \\
1 & x_{1,n}\text{ is not maximal and }x_{2,n}\text{ is minimal} \\
0 & \text{otherwise,} \end{cases} \]
\[ \lambda^{11}_n(x_1,x_2)=\begin{cases}
\lambda^{11}_{n-1}(x_1,x_2) & 
x_{1,n}\text{ is minimal and }x_{2,n}\text{ is minimal} \\
1 & x_{1,n}\text{ is minimal and }x_{2,n}\text{ is not minimal} \\
0 & \text{otherwise} \end{cases} \]
and 
\[ \lambda^{10}_n(x_1,x_2)=\begin{cases}
\lambda^{10}_{n-1}(x_1,x_2) & 
x_{1,n}\text{ is minimal and }x_{2,n}\text{ is maximal} \\
1 & x_{1,n}\text{ is not minimal and }x_{2,n}\text{ is maximal} \\
0 & \text{otherwise.} \end{cases} \]
It is easily checked that $\lambda^{ij}_n$ 
is well-defined and continuous. 

The following is an immediate consequence of 
the definition of $\lambda^{ij}_n$. 

\begin{lem}\label{lambda}
Let $(i,j)\in\{0,1\}^2$. 
For $((x_1,x_2),(y_1,y_2))\in{\cal R}_n$, if 
\[ \lambda^{ij}_n(x_1,x_2)=\lambda^{ij}_n(y_1,y_2), \]
then we have 
\[ \lambda^{ij}_m(x_1,x_2)=\lambda^{ij}_m(y_1,y_2), \]
for all $m>n$. 
\end{lem}

For every $n\in\N$, 
we define a subset ${\cal R}'_n$ of ${\cal R}_n$ by 
\[ {\cal R}'_n=\{((x_1,x_2),(y_1,y_2))\in{\cal R}_n:
\lambda^{ij}_n(x_1,x_2)=\lambda^{ij}_n(y_1,y_2)
\text{ for all }i,j=0,1\}. \]

\begin{lem}
For every $n\in\N$, 
${\cal R}'_n$ is a compact open subequivalence relation 
of ${\cal R}_n$, 
and ${\cal R}'_n$ is contained in ${\cal R}'_{n+1}$. 
\end{lem}
\begin{proof}
It is obvious that ${\cal R}'_n$ is a subequivalence relation 
of ${\cal R}_n$. 
Since $\lambda^{ij}_n$ is continuous, 
${\cal R}'_n$ is compact and open. 
From the lemma above we can see ${\cal R}'_n\subset{\cal R}'_{n+1}$. 
\end{proof}

Define 
\[ {\cal R}'=\bigcup_{n\in\N}{\cal R}'_n. \]
By the lemma above, 
${\cal R}'$ is an AF equivalence relation on $X$. 

\begin{lem}
Let $((x_1,x_2),(y_1,y_2))\in{\cal R}$. 
\begin{enumerate}
\item If $x_1$ is not in $\{\phi_1^n(p_1):n\in\Z\}$, then 
$((x_1,x_2),(y_1,y_2))\in {\cal R}'$. 
\item If $x_2$ is not in $\{\phi_2^n(p_2):n\in\Z\}$, then 
$((x_1,x_2),(y_1,y_2))\in {\cal R}'$. 
\end{enumerate}
\end{lem}
\begin{proof}
It suffices to show (1). 
There exists $n\in\N$ 
such that $((x_1,x_2),(y_1,y_2))\in {\cal R}_n$. 
We can find a natural number $m>n$ 
such that $x_{1,m}$ is not maximal. 
Then, $\lambda^{00}_m(x_1,x_2)$ equals zero. 
From $x_{1,m}=y_{1,m}$, 
we get $\lambda^{00}_m(x_1,x_2)=\lambda^{00}_m(y_1,y_2)=0$. 
It is easy to see that 
$\lambda^{01}_m(x_1,x_2)$ depends only on $x_{2,m}$, and so 
we have $\lambda^{01}_m(x_1,x_2)=\lambda^{01}_m(y_1,y_2)$. 

We can find a natural number $l>n$ 
such that $x_{1,l}$ is not minimal. 
It is clear that $\lambda^{11}_l(x_1,x_2)=0$ and 
$\lambda^{10}_l(x_1,x_2)$ depends only on $x_{2,l}$. 
In a similar fashion to the preceding paragraph, 
we get 
$\lambda^{11}_l(x_1,x_2)=\lambda^{11}_l(y_1,y_2)$ and 
$\lambda^{10}_l(x_1,x_2)=\lambda^{10}_l(y_1,y_2)$. 

By virtue of Lemma \ref{lambda}, 
we can conclude that $((x_1,x_2),(y_1,y_2))$ is 
in ${\cal R}'_k$, 
where $k$ is the maximum of $m$ and $l$. 
\end{proof}

Put $p=(p_1,p_2)\in X$. 
The above lemma tells us that the four ${\cal R}$-orbits 
$[p]_{\cal R}$, $[\phi^{(1,0)}(p)]_{\cal R}$, 
$[\phi^{(0,1)}(p)]_{\cal R}$ and $[\phi^{(1,1)}(p)]_{\cal R}$ 
may split in ${\cal R}'$, 
but the other ${\cal R}$-orbits do not split in ${\cal R}'$. 

\begin{lem}
The equivalence relation ${\cal R}'$ is minimal. 
\end{lem}
\begin{proof}
Let $(x_1,x_2)\in X$. 
It suffices to show that 
$[(x_1,x_2)]_{{\cal R}'}$ is dense in $X$. 
If $x_1$ does not belong to $\{\phi_1^n(p_1):n\in\Z\}$ or 
$x_2$ does not belong to $\{\phi_2^n(p_2):n\in\Z\}$, 
then we have nothing to do, 
because the ${\cal R}'$-orbit of $(x_1,x_2)$ is equal to 
the ${\cal R}$-orbit of it. 
Suppose that $(x_1,x_2)$ is in $\{\phi^a(p):a\in\Z^2\}$. 
Without loss of generality, we may assume that 
$(x_1,x_2)$ belongs to $[p]_{\cal R}$. 
Take finite paths $(e_{1,1},e_{1,2},\dots,e_{1,n})$ in $B_1$ 
and $(e_{2,1},e_{2,2},\dots,e_{2,n})$ in $B_2$. 
Thus $e_{i,k}\in E_{i,k}$ and $r(e_{i,k})=s(e_{i,k+1})$. 
We can find $m>n+2$ such that 
both $x_{1,m}$ and $x_{2,m}$ are maximal. 
It follows that 
$\lambda^{01}_m(x_1,x_2)=0$, $\lambda^{11}_m(x_1,x_2)=0$ 
and $\lambda^{10}_m(x_1,x_2)=1$. 
We have two possibilities: $\lambda^{00}_m(x_1,x_2)=0$ or $1$. 

Let us consider the case that $\lambda^{00}_m(x_1,x_2)$ is one. 
We can find edges $e_{i,k}\in E_{i,k}$ 
for $i=1,2$ and $k=n+1,n+2,\dots,m-1$ such that 
the following are satisfied. 
\begin{itemize}
\item $r(e_{i,k})=s(e_{i,k+1})$ and $r(e_{i,m-1})=s(x_{i,m})$ 
for all $i=1,2$ and $k=n,n+1,\dots,m-2$. 
\item $e_{1,m-1}$ is maximal and $e_{2,m-1}$ is not maximal. 
\end{itemize}
Put 
\[ x'_i=(e_{i,1},e_{i,2},\dots,e_{i,n},e_{i,n+1},\dots,
e_{i,m-1},x_{i,m},x_{i,m+1}\dots)\in X_i \]
for each $i=1,2$. 
Then it is clear that $((x_1,x_2),(x'_1,x'_2))\in {\cal R}_m$. 
Moreover, it is not hard to see 
$\lambda^{00}_m(x'_1,x'_2)=1$, $\lambda^{01}_m(x'_1,x'_2)=0$, 
$\lambda^{11}_m(x'_1,x'_2)=0$ and $\lambda^{10}_m(x'_1,x'_2)=1$. 
Therefore we get $((x_1,x_2),(x'_1,x'_2))\in {\cal R}'_m$. 

Suppose that $\lambda^{00}_m(x_1,x_2)$ is zero. 
In this case we choose the edges $e_{i,k}\in E_{i,k}$ so that 
the following are satisfied. 
\begin{itemize}
\item $r(e_{i,k})=s(e_{i,k+1})$ and $r(e_{i,m-1})=s(x_{i,m})$ 
for all $i=1,2$ and $k=n,n+1,\dots,m-2$. 
\item $e_{1,m-1}$ is not maximal. 
\end{itemize}
Then, we can again obtain $((x_1,x_2),(x'_1,x'_2))\in {\cal R}'_m$. 
Hence we can conclude that 
the ${\cal R}'$-orbit of $(x_1,x_2)$ is dense in $X$. 
\end{proof}

\begin{lem}\label{split}
For every $m\in\N\setminus\{1\}$ we have the following. 
\begin{enumerate}
\item $\lim_{n\to\infty}
(\lambda^{00}_n\times\lambda^{01}_n\times
\lambda^{11}_n\times\lambda^{10}_n)(p_1,\phi_2^{1-m}(p_2))
=(1,0,0,1)$. 
\item $\lim_{n\to\infty}
(\lambda^{00}_n\times\lambda^{01}_n\times
\lambda^{11}_n\times\lambda^{10}_n)(\phi_1(p_1),\phi_2^{1-m}(p_2))
=(0,0,1,0)$. 
\item $\lim_{n\to\infty}
(\lambda^{00}_n\times\lambda^{01}_n\times
\lambda^{11}_n\times\lambda^{10}_n)(p_1,\phi_2^m(p_2))
=(1,0,0,0)$. 
\item $\lim_{n\to\infty}
(\lambda^{00}_n\times\lambda^{01}_n\times
\lambda^{11}_n\times\lambda^{10}_n)(\phi_1(p_1),\phi_2^m(p_2))
=(0,1,1,0)$. 
\item $\lim_{n\to\infty}
(\lambda^{00}_n\times\lambda^{01}_n\times
\lambda^{11}_n\times\lambda^{10}_n)(\phi_1^{1-m}(p_1),p_2)
=(0,0,0,1)$. 
\item $\lim_{n\to\infty}
(\lambda^{00}_n\times\lambda^{01}_n\times
\lambda^{11}_n\times\lambda^{10}_n)(\phi_1^{1-m}(p_1),\phi_2(p_2))
=(1,1,0,0)$. 
\item $\lim_{n\to\infty}
(\lambda^{00}_n\times\lambda^{01}_n\times
\lambda^{11}_n\times\lambda^{10}_n)(\phi_1^m(p_1),p_2)
=(0,0,1,1)$. 
\item $\lim_{n\to\infty}
(\lambda^{00}_n\times\lambda^{01}_n\times
\lambda^{11}_n\times\lambda^{10}_n)(\phi_1^m(p_1),\phi_2(p_2))
=(0,1,0,0)$. 
\end{enumerate}
\end{lem}
\begin{proof}
Straightforward computation. 
\end{proof}

Take a clopen subset $U_i\subset X_i$ which does not contain 
$p_i$ and $\phi_i(p_i)$ for each $i=1,2$. 
Put 
\[ B=(\{p_1\}\times U_2)\cup(U_1\times\{p_2\}) \]
and 
\[ B^*=(\{\phi_1(p_1)\}\times U_2)\cup(U_1\times\{\phi_2(p_2)\}). \]

\begin{lem}\label{etale}
Both $B$ and $B^*$ are closed ${\cal R}'$-\'etale thin subsets. 
\end{lem}
\begin{proof}
It suffices to show the statement for $B$. 
Suppose that $((x_1,x_2),(y_1,y_2))$ is in ${\cal R}'\cap(B\times B)$. 
Without loss of generality, we may assume $x_1=p_1$. 
Suppose $y_2=p_2$. 
Then $x_2$ must be $\phi_2^{1-m}(p_2)$ for some $m\in\N$, 
and $m$ is not one because $p_2$ is not in $U_2$. 
Similarly $y_1$ must be $\phi_1^{1-l}(p_1)$ 
for some $l\in\N\setminus\{1\}$. 
But $((p_1,\phi_1^{1-m}(p_2)),(\phi_1^{1-l}(p_1),p_2))$ never be 
in ${\cal R}'$ by the lemma above. 
Hence we have $y_1=p_1$. 
Thus $((x_1,x_2),(y_1,y_2))$ is equal to 
$((p_1,x_2),(p_1,\phi_2^m(x_2)))$ for some $m\in\Z$. 
Define 
\[ V=\{((a,b),(c,d))\in {\cal R}': a=c, \ d=\phi_2^m(b)
\text{ and }b,d\in U_2\}. \]
Then $V$ is a clopen neighborhood of 
$((p_1,x_2),(p_1,\phi_2^m(x_2)))$ in ${\cal R}'$. 
For $((a,b),(c,d))\in V$, it is obvious that 
$(a,b)\in B$ if and only if $(c,d)\in B$, 
which implies that $B$ is \'etale. 

We would like to show that a probability measure on $X=X_1\times X_2$ is 
${\cal R}$-invariant if and only if it is ${\cal R}'$-invariant. 
If this is shown, thinness of $B$ easily follows. 
But, except for countably many $(x_1,x_2)$'s, 
the equivalence class $[(x_1,x_2)]_{\cal R}$ is equal to 
$[(x_1,x_2)]_{{\cal R}'}$. 
Since every invariant measure is nonatomic, we can finish the proof. 
\end{proof}

\begin{lem}
We have ${\cal R}'\cap (B\times B^*)=\emptyset$. 
\end{lem}
\begin{proof}
Suppose that $((x_1,x_2),(y_1,y_2))$ is contained 
in ${\cal R}'\cap (B\times B^*)$. 
Without loss of generality, we may assume $x_1=p_1$. 
Then $y_1$ never be $\phi_1(x_1)$, 
because $((p_1,x_2),(\phi_1(p_1),y_2))$ does not belong to ${\cal R}$. 
It follows that $y_2=\phi_2(p_2)$ and 
$((x_1,x_2),(y_1,y_2))=
((p_1,\phi_2^m(p_2)),(\phi_1^{1-l}(p_1),\phi_2(p_2)))$ 
for some $m,l\in\N\setminus\{1\}$. 
This pair, however, never belongs to ${\cal R}'$ 
by virtue of Lemma \ref{split}, which completes the proof. 
\end{proof}

We define a homeomorphism $\beta:B\to B^*$ as follows. 
For $(p_1,x_2)\in \{p_1\}\times U_2$, 
we put $\beta(p_1,x_2)=(\phi_1(p_1),x_2)$. 
For $(x_1,p_2)\in U_1\times\{p_2\}$, 
we put $\beta(x_1,p_2)=(x_1,\phi_2(p_2))$. 

\begin{lem}
The homeomorphism $\beta:B\to B^*$ induces an isomorphism 
between ${\cal R}'\cap(B\times B)$ and ${\cal R}'\cap(B^*\times B^*)$. 
\end{lem}
\begin{proof}
Since the topology of 
${\cal R}'\cap(B\times B)$ and ${\cal R}'\cap(B^*\times B^*)$ 
is inherited from ${\cal R}$, 
it suffices to show that $\beta$ is a well-defined bijection 
between ${\cal R}'\cap(B\times B)$ and ${\cal R}'\cap(B^*\times B^*)$. 
Let $((x_1,x_2),(y_1,y_2))\in{\cal R}'\cap(B\times B)$. 
Without loss of generality, 
we may assume that $x_1=p_1$ and $x_2\in U_2$. 
By the proof of Lemma \ref{etale}, 
we get $y_1=p_1$ and $y_2\in U_2$. 
It follows that $\beta(p_1,x_2)=(\phi_1(p_1),x_2)$ and 
$\beta(p_1,y_2)=(\phi_1(p_1),y_2)$. 
If $x_2$ does not belong to $\{\phi_2^n(p_2):n\in\Z\}$, 
then the ${\cal R}'$-orbit of $(\phi_1(p_1),x_2)$ is equal to 
the ${\cal R}$-orbit of it. 
Hence we have 
$((\phi_1(p_1),x_2),(\phi_1(p_1),y_2))\in{\cal R}'\cap(B^*\times B^*)$. 

Suppose that $x_2$ and $y_2$ belong to $\{\phi_2^n(p_2):n\in\Z\}$. 
Since $((p,x_2),(p,y_2))\in{\cal R}'\subset{\cal R}$, 
we have two possibilities: 
both $x_2$ and $y_2$ belong to $\{\phi_2^{1-n}(p_2):n\in\N\}$, or 
both $x_2$ and $y_2$ belong to $\{\phi_2^n(p_2):n\in\N\}$. 
Without loss of generality, we may assume the latter. 
Thus, $x_2=\phi_2^n(p_2)$ and 
$y_2=\phi_2^m(p_2)$ for some $n,m\in\N$. 
Because $x_2$ and $y_2$ is in $U_2$, $n$ and $m$ are greater than one. 
It follows from Lemma \ref{split} that 
$((\phi_1(p_1),\phi_2^n(p_2)),(\phi_1(p_1),\phi_2^m(p_2)))$ 
belongs to ${\cal R}'$. 
The proof is completed. 
\end{proof}

\begin{lem}
Let $\tilde{\cal R}$ be the equivalence relation generated 
by ${\cal R}'$ and the graph of $\beta$. 
Then ${\cal R}_\phi$ is generated by $\tilde{\cal R}$ and 
$(p,\phi^{(0,1)}(p))$, $(\phi^{(0,1)}(p),\phi^{(1,1)}(p))$ 
and $(\phi^{(1,1)}(p),\phi^{(1,0)}(p))$. 
\end{lem}
\begin{proof}
Evidently ${\cal R}_\phi$ is generated 
by ${\cal R}$ and the graph of $\beta$. 
As mentioned before, 
if $(x_1,x_2)$ is not contained 
in the ${\cal R}_\phi$-orbit of $p=(p_1,p_2)$, 
then its ${\cal R}$-orbit agrees with its ${\cal R}'$-orbit. 
It follows that the $\tilde{\cal R}$-orbit of $(x_1,x_2)$ 
agrees with the ${\cal R}_\phi$-orbit of it. 

Let us consider $[p]_{{\cal R}_\phi}$. 
Notice that it splits into four orbits in ${\cal R}$, 
namely the ${\cal R}$-orbits of 
$p$, $\phi^{(0,1)}(p)$, $\phi^{(1,0)}(p)$ and $\phi^{(1,1)}(p)$. 
From Lemma \ref{split} we can see that these orbits split 
into eight orbits in ${\cal R}'$, namely 
the ${\cal R}'$-orbits of 
$p$, $\phi^{(-1,0)}(p)$, $\phi^{(0,1)}(p)$, $\phi^{(0,2)}(p)$, 
$\phi^{(1,0)}(p)$, $\phi^{(1,-1)}(p)$, 
$\phi^{(1,1)}(p)$ and $\phi^{(2,1)}(p)$. 
It can be easily seen that 
\[ [p]_{\tilde{\cal R}}
=[p]_{{\cal R}'}\cup[\phi^{(1,-1)}(p)]_{{\cal R}'}, \]
\[ [\phi^{(0,1)}(p)]_{\tilde{\cal R}}
=[\phi^{(0,1)}(p)]_{{\cal R}'}\cup[\phi^{(-1,0)}(p)]_{{\cal R}'}, \]
\[ [\phi^{(1,0)}(p)]_{\tilde{\cal R}}
=[\phi^{(1,0)}(p)]_{{\cal R}'}\cup[\phi^{(2,1)}(p)]_{{\cal R}'} \]
and 
\[ [\phi^{(1,1)}(p)]_{\tilde{\cal R}}
=[\phi^{(1,1)}(p)]_{{\cal R}'}\cup[\phi^{(0,2)}(p)]_{{\cal R}'}. \]
Therefore, by glueing the $\tilde{\cal R}$-orbits of 
$p$, $\phi^{(0,1)}(p)$, $\phi^{(1,0)}(p)$ and $\phi^{(1,1)}(p)$, 
we can recover the equivalence relation ${\cal R}_\phi$. 
\end{proof}

\bigskip

By \cite[Theorem 4.6]{HPS}, 
every minimal homeomorphism on the Cantor set is conjugate to 
a Bratteli-Vershik transformation 
on a simple properly ordered Bratteli diagram. 
Hence we can summarize 
the results obtained in this section as follows. 

\begin{thm}\label{summary}
Let $(X_1,\phi_1)$ and $(X_2,\phi_2)$ be two Cantor minimal $\Z$-systems 
and let $p_1\in X_1$ and $p_2\in X_2$. 
Take clopen subsets $U_1\subset X_1$ and $U_2\subset X_2$ 
so that $p_i$ and $\phi_i(p_i)$ do not belong to $U_i$ for each $i=1,2$. 
Put $B=(\{p_1\}\times U_2)\cup(U_1\times\{p_2\})$ and 
$B^*=(\{\phi_1(p_1)\}\times U_2)\cup(U_1\times\{\phi_2(p_2)\})$. 
Define $\beta:B\to B^*$ 
by $\beta(p_1,x_2)=(\phi_1(p_1),x_2)$ and 
$\beta(x_1,p_2)=(x_1,\phi_2(p_2))$. 
Let $\phi$ be the $\Z^2$-action on $X=X_1\times X_2$ 
induced by $\phi_1\times\id$ and $\id\times\phi_2$. 
Put $p=(p_1,p_2)$. 

Then we can find a subequivalence relation 
${\cal R}'\subset{\cal R}_\phi$ 
such that the following are satisfied. 
\begin{enumerate}
\item ${\cal R}'$ is a minimal AF equivalence relation, 
where the topology is given by ${\cal R}_\phi$. 
\item Both $B$ and $B^*$ are closed ${\cal R}'$-\'etale thin subsets. 
\item ${\cal R}'\cap(B\times B^*)$ is empty. 
\item $\beta:B\to B^*$ induces an isomorphism 
between ${\cal R}'\cap(B\times B)$ and ${\cal R}'\cap(B^*\times B^*)$. 
\item The equivalence relation ${\cal R}_\phi$ is 
generated by ${\cal R}'$, the graph of $\beta$ and 
\[ \{(p,\phi^{(0,1)}(p)), \ (\phi^{(0,1)}(p),\phi^{(1,1)}(p)) \ 
(\phi^{(1,1)}(p),\phi^{(1,0)}(p))\}. \]
\end{enumerate}
\end{thm}

\section{The main result}

Let $\xp$ and $\yp$ be two Cantor minimal $\Z^2$-systems and 
let $\pi:\yp\to\xp$ be a factor map. 

\begin{lem}\label{extAF}
Suppose that ${\cal R}$ is an open subequivalence relation of 
${\cal R}_\phi$. 
For 
\[ {\cal S}=\{(y,y')\in{\cal R}_\psi:(\pi(y),\pi(y'))\in{\cal R}\}, \]
we have the following. 
\begin{enumerate}
\item If ${\cal R}$ is compact and open, 
then ${\cal S}$ is also compact and open. 
\item If ${\cal R}$ is AF, then ${\cal S}$ is also AF.
\end{enumerate}
\end{lem}
\begin{proof}
(2) follows immediately from (1). 
Suppose that ${\cal R}$ is compact and open. 
Since $\pi\times\pi:{\cal R}_\psi\to{\cal R}_\phi$ is 
proper and continuous, 
we can see that ${\cal S}={\cal R}_\psi\cap(\pi\times\pi)^{-1}({\cal R})$ 
is compact and open. 
\end{proof}

\begin{lem}\label{extetale}
Suppose that ${\cal R}$ is an open subequivalence relation of 
${\cal R}_\phi$. 
Let ${\cal S}={\cal R}_\psi\cap(\pi\times\pi)^{-1}({\cal R})$. 
If $B\subset X$ is a closed ${\cal R}$-\'etale thin subset, 
then $\pi^{-1}(B)$ is a closed ${\cal S}$-\'etale thin subset. 
\end{lem}
\begin{proof}
Let $\mu$ be an ${\cal S}$-invariant probability measure. 
Then we have 
\[ \mu(\pi^{-1}(B))=\pi_*(\mu)(B)=0, \]
because $\pi_*(\mu)$ is a ${\cal R}$-invariant probability measure. 

Take $(y,y')\in {\cal S}$. 
By the \'etaleness of $B$, 
we can find a clopen neighborhood $V$ of $(\pi(y),\pi(y'))$ 
in ${\cal R}$ such that, for $(x,x')\in V$, 
we have $x\in B$ if and only if $x'\in B$. 
It is clear that $U={\cal R}_\psi\cap(\pi\times\pi)^{-1}(V)$ is 
a clopen neighborhood of $(y,y')$ in ${\cal S}$. 
Suppose $(z,z')\in U$. 
Because of $(\pi(z),\pi(z'))\in V$, we have 
\[ z\in \pi^{-1}(B)\Leftrightarrow\pi(z)\in B
\Leftrightarrow\pi(z')\in B\Leftrightarrow z'\in\pi^{-1}(B). \]
It follows that $\pi^{-1}(B)$ is ${\cal S}$-\'etale. 
\end{proof}

Now we are ready to prove the main theorem. 

\begin{thm}
Let $\pi:\yp\to\xp$ be a factor map 
between Cantor minimal $\Z^2$-systems. 
If $\xp$ is conjugate to a product of two Cantor minimal $\Z$-systems, 
then ${\cal R}_\psi$ is affable. 
\end{thm}
\begin{proof}
We may assume that $\xp$ is equal to 
the product of two Cantor minimal $\Z$-systems 
$(X_1,\phi_1)$ and $(X_2,\phi_2)$, that is, 
$X=X_1\times X_2$ and 
$\phi^{(n,m)}(x_1,x_2)=(\phi_1^n(x_1),\phi_2^m(x_2))$ 
for all $(n,m)\in\Z^2$.
Let $p=(p_1,p_2)$, $U_1,U_2$, $B$, $B^*$, 
$\beta:B\to B^*$ and ${\cal R}'$ 
be as in Theorem \ref{summary}. 

Put ${\cal S}={\cal R}_\psi\cap(\pi\times\pi)^{-1}({\cal R}')$. 
Thanks to Theorem \ref{summary} (1) and Lemma \ref{extAF}, 
the equivalence relation ${\cal S}$ is AF, 
where the topology is given by ${\cal R}_\psi$. 
In order to show that ${\cal S}$ is minimal, 
let us choose $x_i\in X_i\setminus\{p_i\}$ and 
put $x_0=(x_1,x_1)\in X$. 
Take $y\in Y$ arbitrarily. 
The closure of $[\pi(y)]_{{\cal R}'}$ is $X$, 
because ${\cal R}'$ is minimal. 
It follows that the closure of $[y]_{\cal S}$ contains 
a preimage of $x_0$, namely $y_0\in Y$. 
On account of $[y_0]_{{\cal R}_\psi}=[y_0]_{\cal S}$, 
we can see that $[y_0]_{\cal S}$ is dense in $Y$. 
Therefore $[y]_{\cal S}$ is dense in $Y$. 

Put $C=\pi^{-1}(B)$ and $C^*=\pi^{-1}(B^*)$. 
By means of Theorem \ref{summary} (2) and Lemma \ref{extetale}, 
we have that 
both $C$ and $C^*$ are closed ${\cal S}$-\'etale thin subsets. 
Moreover, 
it is easily seen that ${\cal S}\cap(C\times C^*)$ is empty. 

We define a homeomorphism $\gamma:C\to C^*$ as follows. 
Take $y\in C$. 
If $\pi(y)=(p_1,x_2)$ for some $x_2\in U_2$, 
then we set $\gamma(y)=\psi^{(1,0)}(y)$. 
If $\pi(y)=(x_1,p_2)$ for some $x_1\in U_1$, 
then we set $\gamma(y)=\psi^{(0,1)}(y)$. 
It is routine to check that $\gamma$ is a well-defined homeomorphism 
from $C$ to $C^*$ and $\gamma$ induces an isomorphism 
between ${\cal S}\cap(C\times C)$ and ${\cal S}\cap(C^*\times C^*)$. 

Let $\tilde{\cal S}$ be the equivalence relation 
generated by ${\cal S}$ and the graph of $\gamma$. 
We can apply the absorption theorem \cite[Theorem 4.18]{GPS2} 
to ${\cal S}$ and $\gamma:C\to C^*$ 
and get that $\tilde{\cal S}$ is affable. 

The equivalence relation $\tilde{\cal S}$ is a little smaller 
than ${\cal R}_\psi$. 
We resolve this problem 
by using the absorption theorem three more times. 
Let $D_1=\pi^{-1}(p)$, $D_2=\pi^{-1}(\phi^{(0,1)}(p))$, 
$D_3=\pi^{-1}(\phi^{(1,1)}(p))$ and 
$D_4=\pi^{-1}(\phi^{(1,0)}(p))$. 
At first, 
we apply the absorption theorem to $\psi^{(0,1)}:D_1\to D_2$. 
Notice that 
\[ \tilde{\cal S}\cap(D_i\times D_i)=\{(y,y):y\in D_i\} \]
for each $i=1,2$ and that $\tilde{\cal S}\cap(D_1\times D_2)$ is empty. 
Therefore the hypothesis of the absorption theorem is 
trivially satisfied. 
It follows that the equivalence relation 
generated by $\tilde{\cal S}$ and 
\[ \{(y,\psi^{(0,1)}(y)):y\in D_1\} \]
is affable. 
Theorem \ref{summary} (5) and 
two more applications of the absorption theorem imply 
that ${\cal R}_\psi$ is affable. 
\end{proof}

\flushleft{
\textit{e-mail: matui@math.s.chiba-u.ac.jp \\
Graduate School of Science and Technology,\\
Chiba University,\\
1-33 Yayoi-cho, Inage-ku,\\
Chiba 263-8522,\\
Japan. }}

\end{document}